\documentclass[12pt,reqno,a4wide]{amsart}

\usepackage{tikz} 
\usepackage{calrsfs}
\usepackage{mathrsfs}

\usepackage{array}
\usepackage{hyperref}

\usetikzlibrary{decorations.pathreplacing}
\usetikzlibrary{fit}

\usepackage{pgfplots}

\allowdisplaybreaks

\catcode`,\active

\catcode`\,12

\begin{document}
\bibliographystyle{plain}

%
%

	\title
	{Grand Motzkin paths and $\{0,1,2\}$-trees -- a simple bijection}

	\author[H. Prodinger ]{Helmut Prodinger }
	\address{Department of Mathematics, University of Stellenbosch 7602, Stellenbosch, South Africa
	and
NITheCS (National Institute for
Theoretical and Computational Sciences), South Africa.}
	\email{hproding@sun.ac.za}

	
	\begin{abstract}
		A well-known bijection between Motzkin paths and ordered trees with outdegree always $\le2$, is lifted to Grand Motzkin paths (the nonnegativity is dropped) 
		and an ordered list of an odd number of such $\{0,1,2\}$ trees. This offers an alternative to a recent paper by Rocha and Pereira Spreafico.
	\end{abstract}
	

\maketitle

\section{Introduction}

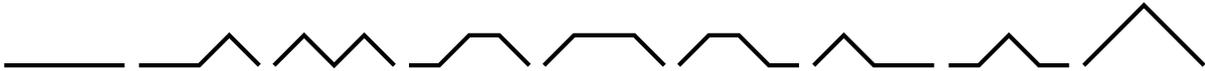
\begin{figure}[h]
	\label{all1}
\begin{center}

			\begin{tikzpicture}[scale=0.4]
				
				\draw[ultra thick] (0,0) to (1,0) to (2,0) to (3,0) to (4,0)  ;
			\end{tikzpicture}
			\begin{tikzpicture}[scale=0.4]
				
				\draw[ ultra thick] (0,0) to (1,0) to (2,0) to (3,1) to (4,0)  ;
			\end{tikzpicture}
			\begin{tikzpicture}[scale=0.4]
				
				\draw[ ultra thick] (0,0) to (1,1) to (2,0) to (3,1) to (4,0)  ;
			\end{tikzpicture}
			\begin{tikzpicture}[scale=0.4]
				
				\draw[ ultra thick] (0,0) to (1,0) to (2,1) to (3,1) to (4,0)  ;
			\end{tikzpicture}
			\begin{tikzpicture}[scale=0.4]
				
				\draw[ ultra thick] (0,0) to (1,1) to (2,1) to (3,1) to (4,0)  ;
			\end{tikzpicture}
			\begin{tikzpicture}[scale=0.4]
				
				\draw[ ultra thick] (0,0) to (1,1) to (2,1) to (3,0) to (4,0)  ;
			\end{tikzpicture}
			\begin{tikzpicture}[scale=0.4]
				
				\draw[ ultra thick] (0,0) to (1,1) to (2,0) to (3,0) to (4,0)  ;
			\end{tikzpicture}
			\begin{tikzpicture}[scale=0.4]
				
				\draw[ ultra thick] (0,0) to (1,0) to (2,1) to (3,0) to (4,0)  ;
			\end{tikzpicture}
			\begin{tikzpicture}[scale=0.4]
				
				\draw[ ultra thick] (0,0) to (1,1) to (2,2) to (3,1) to (4,0)  ;
			\end{tikzpicture}
	 	
\end{center}

\caption{All 9 Motzkin of 4 steps (length 4).}
\end{figure}

Motzkin paths  appear first in \cite{Motzkin48}. In the encyclopedia \cite{OEIS} they 
are enumerated by sequence A001006, with many references given. They consist of up-steps $U=(1,1)$, down-steps $D=(1,-1)$ and horizontal (flat)
steps $F=(1,0)$.  They start at the origin and must never go below the $x$-axis. Usually one
requires the path to end on the $x$-axis as well, but occasionally one  uses the term \emph{Motzkin path} also for paths that end
on a different level.  Figure~\ref{all1} shows all Motzkin paths of 4 steps (=length 4).

The enumeration of Motzkin paths is done using the generating function $M=M(z)$ and a decomposition according to the first return to the $x$-axis, viz.
\begin{equation*}
M=1+zM+z^2M^2,
\end{equation*}
this can be found in many books, e.g. in \cite{FS}. Solving,
\begin{equation*}
M(z)=\frac{1-z-\sqrt{1-2z-3z^2}}{2z^2}.
\end{equation*}

The other combinatorial structure that plays a role in this note are \emph{ordered trees}. They are enumerated by an equation for the generating function (according to the number of nodes)
\begin{equation*}
P=z+zP+zP^2+zP^3+\cdots=\frac{z}{1-P}\quad\text{and therefore}\quad P=P(z)=\frac{1-\sqrt{1-4z}}{2}.
\end{equation*}
The subclass of $\{0,1,2\}$-trees of interest only allows outdegrees $0,1,2$ and so we get again a generating function
\begin{equation*}
	Q=z+zQ+zQ^2\quad\text{and therefore}\quad Q=Q(z)=\frac{1-z-\sqrt{1-2z-3z^2}}{2z}=zM(z).
\end{equation*}
So there should be a bijection of $\{0,1,2\}$-trees with $n+1$ nodes (= $n$ edges) and Motzkin paths of length $n$; the simplest I know is from \cite{DeutschShapiro}:

One runs through the $\{0,1,2\}$-tree in \emph{pre-order}; if one sees an edge for the \emph{first} time, one translates a single edge (degree 1) into a flat step, 
a left edge into an up-step and a right edge into a down-step. It is easy to see that the process is reversible, which is the desired bijection.

\section{Grand Motzkin paths}

As a first step, we need the generating function of Motzkin paths, ending on level $k$, not just the usual case $0$. This is a standard argument, by decomposing such a path according to the last time you visit level $0$, then an up-step, and we wait until we visit level $1$ for the last time, and so one.
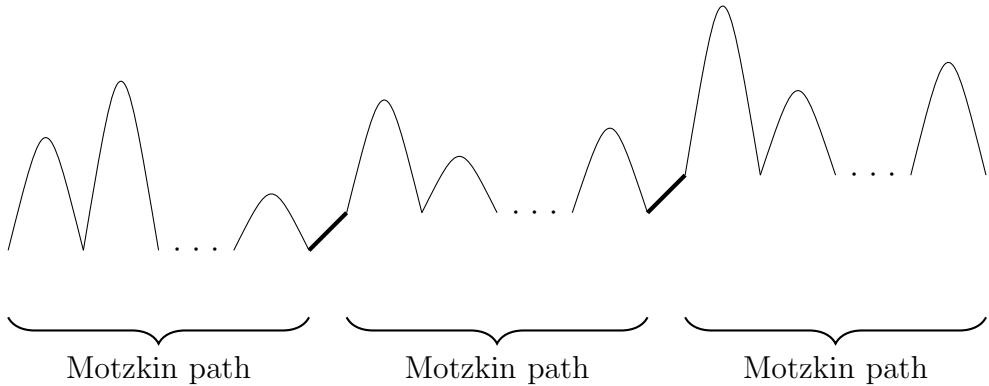
\begin{figure}[h]
	\begin{center}
		\begin{tikzpicture}[scale=0.5]
			
			
			\draw (1,1) .. controls (2,5) .. (3,1); \draw
			(3,1) .. controls (4,7) .. (5,1); \node at
			(5.5,1) {$\cdot$}; \node at (6,1) {$\cdot$};
			\node at (6.5,1) {$\cdot$}; \draw (7,1) ..
			controls (8,3) .. (9,1); \draw [ultra
			thick](9,1) to  (10,2); \draw (1+9,1+1) ..
			controls (2+9,6) .. (3+9,1+1); \draw (3+9,1+1)
			.. controls (4+9,4) .. (5+9,1+1); \node at
			(5.5+9,2) {$\cdot$}; \node at (6+9,2)
			{$\cdot$}; \node at (6.5+9,2) {$\cdot$}; \draw
			(16,2) .. controls (17,5) .. (18,2); 
			\draw 			[ultra thick](18,2) to  (19,3);
			
			\draw (19,0+3) .. controls (20,6+3) .. (21,0+3);
			\draw (21,0+3) .. controls (22,3+3) .. (23,0+3);
			\node at (23.5,0+3) {$\cdot$}; \node at (24,0+3)
			{$\cdot$}; \node at (24.5,0+3) {$\cdot$}; \draw
			(25,0+3) .. controls (26,4+3) .. (27,0+3);

			\draw [thick, decorate, decoration={brace,
				amplitude=10pt, mirror, raise=4pt}] (1cm, -0.5)
			to node[below,yshift=-0.5cm] {Motzkin path} (9cm,
			-0.5);
			
			\draw [thick, decorate, decoration={brace,
				amplitude=10pt, mirror, raise=4pt}] (10cm,
			-0.5) to node[below,yshift=-0.5cm] {Motzkin path}
			(18cm,- 0.5);
			
			\draw [thick, decorate, decoration={brace,
				amplitude=10pt, mirror, raise=4pt}] (19cm,
			-0.5) to node[below,yshift=-0.5cm] {Motzkin path}
			(27cm, -0.5);
			
		\end{tikzpicture} \end    {center}
		
		\caption{The decomposition of a Motzkin path that ends at level $2$.}
	\end{figure}

From this we find the generating function as $z^kM^{k+1}$. This is well-known.

Now we come to Grand Motzkin paths, a notation I picked up from \cite{luca}. This more general family has the same steps as Motzkin paths, returns to the
$x$-axis at the end, but the condition that the paths must stay (weakly) above the $x$-axis is dropped. 

Let $k$ be the unique negative number such that the Grand Motzkin paths reaches the level $-k$ when going from left to right; for $k=0$ they are just ordinary Motzkin paths. We consider the first such point $(a,-k)$ and the last such point $(b,-k)$; it is possible that $a=b$. Then the paths decomposes canonically into 3 parts:
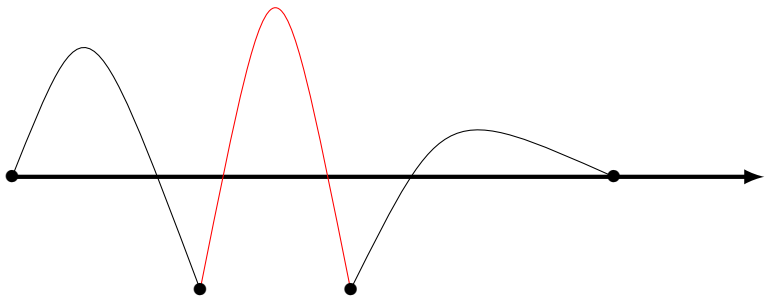
\begin{figure}[h]
	\begin{center}
		\begin{tikzpicture}[scale=0.5]
			\draw[ultra thick,-latex](0,0) to (20,0);
			
			
			\draw (0,0) .. controls (2,5)  .. (5,-3); 
			\draw[red] (5,-3) .. controls (7,7)  .. (9,-3); 
				\draw (9,-3) .. controls (11.5,2)  .. (16,0); 
				
				 \node at (0, 0)  (c2)    {$\bullet$};
				 \node at (16, 0)  (c3)    {$\bullet$};
				 \node at (5, -3)  (c4)    {$\bullet$};
 				 \node at (9, -3)  (c5)    {$\bullet$};
			
		\end{tikzpicture} 
	\end{center}
		
		\caption{The decomposition of a Grand Motzkin path that has $-3$ as its minimal level.}
	\end{figure}
This decomposition produces the generating function of such paths as a product of 3 terms: The first part corresponds   to
$z\cdot z^{k-1}M^{k}$ (the last step must be a down-step), the second part to $M$, and the third part again to $z\cdot z^{k-1}M^{k}$ (the first step must be an up-step). The product is then $z^{2k}M^{2k+1}$. Summing over all possible values of $k$, the enumeration of Grand Motzkin paths is done via the generating function
\begin{equation*}
M+\sum_{k\ge1}{z^{2k}M^{2k+1}}=\sum_{k\ge0}{z^{2k}M^{2k+1}}=\frac1{\sqrt{1-2z-3z^2}}.
\end{equation*} 

In the world of $\{0,1,2\}$-trees, we form a new super-root, with $2k+1$ successors, each of which is a $\{0,1,2\}$-tree.
The generating function is then
\begin{equation*}
z\sum_{k\ge0}Q^{2k+1}=z^2\sum_{k\ge0}z^{2k}M^{2k+1}=\frac{z^2}{\sqrt{1-2z-3z^2}}.
\end{equation*}
The previous bijection takes over to the new situation: if a grand Motzkin path consists of $2k+1$ ordinary Motzkin paths, then each of them corresponds to a 
$\{0,1,2\}$-tree as before. The extra factor $z^2$ stems from the fact that for $\{0,1,2\}$-trees, the edges correspond to the steps. And now, there is the super-root, which should not be counted when considering the corresponding Grand Motzkin path.

The paper \cite{rocha} has a correspondence between Grand Motzkin paths and $\{0,1,2\}$-trees with a super-root with an odd number of successors. However the arguments are perhaps less direct than the present ones.

\section{Trinomial coefficients and enumeration}

The trinomial coefficients (notation from Comtet \cite{Comtet-book}) are given by
\begin{equation*}
\binom{n,3}{k}=[z^k](1+z+z^2)^n.
\end{equation*}
These coefficients are intimately related to the Motzkin-world, as we will discuss for the reader's benefit. We use the substitution $z=\dfrac{v}{1+v+v^2}$ as we first did in \cite{Prodinger-three}. Using this substitution, all generating functions become much easier, like 
\begin{equation*}
Q(z)=v,\quad \frac1{\sqrt{1-2z-3z^2}}=\frac{1+v+v^2}{1-v^2}.
\end{equation*}
Coefficients can be extracted via contour integration, which is a variant of the Lagrange inversion formula. See \cite{BrKnRi72} and \cite{Prodinger-three} as illustrations of the technique. Note that when $z$ runs around the origin in a small circle, $v$ runs around the origin once as well, in a deformed circle. We will show two sample computations:
\begin{align*}
[z^n]\frac1{\sqrt{1-2z-3z^2}}&=\frac1{2\pi i}\oint \frac{dz}{z^{n+1}}\frac1{\sqrt{1-2z-3z^2}}\\
&=\frac1{2\pi i}\oint \frac{dv(1-v^2)}{(1+v+v^2)^2}\frac{(1+v+v^2)^{n+1}}{v^{n+1}}\frac{1+v+v^2}{1-v^2}\\
&=\frac1{2\pi i}\oint \frac{dv}{v^{n+1}}(1+v+v^2)^{n}=[v^n](1+v+v^2)^{n}=\binom{n,3}{n};
\end{align*}
and
\begin{align*}
[z^n]Q^j&=\frac1{2\pi i}\oint \frac{dz}{z^{n+1}} v^j=\frac1{2\pi i}\oint \frac{dv(1-v^2)}{(1+v+v^2)^2}\frac{(1+v+v^2)^{n+1}}{v^{n+1}} v^j\\
&=\frac1{2\pi i}\oint \frac{dv(1-v^2)}{v^{n+1-j}} (1+v+v^2)^{n-1} \\
&=[v^{n-j}](1+v+v^2)^{n-1}-[v^{n-j-2}](1+v+v^2)^{n-1}=\binom{n-1,3}{n-j}-\binom{n-1,3}{n-j-2}.
\end{align*}

\bibliographystyle{plain}

\end{document}